\documentclass[11pt]{article}

\usepackage{fullpage}
\usepackage{xypic}
\usepackage{amsgen}
\usepackage{amsmath}
\usepackage{amstext}
\usepackage{amsbsy}
\usepackage{amsopn}
\usepackage{amsfonts}
\usepackage{amssymb}
\usepackage{eepic}
\usepackage{graphicx}
\usepackage{epsf}
\usepackage{pstricks}
\xyoption{all}

\def\Box{\square}

\def\mapright#1{\smash{\mathop{\longrightarrow}\limits^{#1}}}
\def\tra#1{\smash{\mathop{\mid\kern
-1pt\joinrel\relbar\joinrel\relbar}\limits^{*}_{#1}}}
\def\longtra#1{\smash{\mathop{\mid\kern
-1pt\joinrel\relbar\joinrel\relbar\joinrel\relbar}\limits^{*}_{#1}}}
\def\vlongtra#1{\smash{\mathop{\mid\kern
-1pt\joinrel\relbar\joinrel\relbar\joinrel\relbar\joinrel\relbar}\limits^{*}_{#1}}}
\def\vvlongtra#1{\smash{\mathop{\mid\kern
-1pt\joinrel\relbar\joinrel\relbar\joinrel\relbar\joinrel\relbar\joinrel\relbar}\limits^{*}_{#1}}}
\def\vvvlongtra#1{\smash{\mathop{\mid\kern
-1pt\joinrel\relbar\joinrel\relbar\joinrel\relbar\joinrel\relbar\joinrel\relbar\joinrel\relbar}\limits^{*}_{#1}}}
\def\etra#1{\smash{\mathop{\mid\kern
-1pt\joinrel\relbar\joinrel\relbar}\limits_{#1}}}
\def\mapleft#1{\smash{\mathop{\longleftarrow}\limits^{#1}}}

\def\A{{\cal{A}}}
\def\iff{\Leftrightarrow}
\def\Rw{\Rightarrow}

\def\C{{\cal{C}}}
\def\F{{\cal{F}}}

\def\N{\mathbb{N}}

\def\S{{\cal{S}}}

\def\fix{\mbox{Fix}}
\def\per{\mbox{Per}}

\def\ker{\mbox{Ker}\,}

\def\endo{\mbox{End}}

\def\min{\mbox{min}}
\def\rk{\mbox{rk}}
\def\krk{\mbox{Krk}}
\def\prk{\mbox{prk}}

\def\sup{\mbox{sup}}

\def\G{{\cal{G}}}

\def\Z{\mathbb{Z}}

\def\p{\varphi}

\def\inv{^{-1}}

\def\bi{\begin{itemize}}
\def\ei{\end{itemize}}
\def\beq{\begin{equation}}
\def\eeq{\end{equation}}

\newtheorem{T}{Theorem}[section]
\newcommand{\bt}{\begin{T}}
\newcommand{\et}{\end{T}}
\newcommand{\ftd}{$\square$\end{T}}

\newtheorem{Proposition}[T]{Proposition}
\newcommand{\bp}{\begin{Proposition}}
\newcommand{\ep}{\end{Proposition}}
\newcommand{\fpd}{$\square$\end{Proposition}}

\newtheorem{Lemma}[T]{Lemma}
\newcommand{\bl}{\begin{Lemma}}
\newcommand{\el}{\end{Lemma}}
\newcommand{\fld}{$\square$\end{Lemma}}

\newtheorem{Corol}[T]{Corollary}
\newcommand{\bc}{\begin{Corol}}
\newcommand{\ec}{\end{Corol}}
\newcommand{\fcd}{$\square$\end{Corol}}

\newtheorem{Result}[T]{Result}
\newcommand{\br}{\begin{Result}}
\newcommand{\er}{\end{Result}}
\newcommand{\frd}{$\square$\end{Result}}

\newtheorem{Example}[T]{Example}
\newcommand{\be}{\begin{Example}}
\newcommand{\ee}{\end{Example}}

\newtheorem{Problem}[T]{Problem}
\newcommand{\bq}{\begin{Problem}}
\newcommand{\eq}{\end{Problem}}

\newcommand{\proof}
   {\par\medbreak\noindent{\bf Proof}.\enspace}

\newcommand{\qed}{
$\Box$
\par\bigbreak}

\normalbaselines
\def\abstract#1{\par\bigskip
\begingroup\small
\baselineskip=12truept
\begin{center}ABSTRACT\end{center}
\par\medskip\par\noindent
\null\hfill\hbox{\vbox{\hsize=5truein\noindent#1}}
\hfill\null\par\endgroup\par}

\title{Finiteness results for subgroups of finite extensions}
\author{{\bf V\'\i tor Ara\'ujo}\\
$ $\\ {\em Universidade Federal da Bahia, Instituto de
  Matem\'atica,}\\
{\em Av. Adhemar de Barros, S/N, Ondina,}\\
{\em 40170-110 Salvador-BA, Brazil}\\
{\em email:} vitor.d.araujo@ufba.br\\
$ $\\
{\bf Pedro V. Silva}\\ $ $\\
{\em Universidade Federal da Bahia, Instituto de
  Matem\'atica, Brazil}\\
{\em and Centro de
Matem\'{a}tica, Faculdade de Ci\^{e}ncias, Universidade do
Porto,}\\ {\em R. Campo Alegre 687, 4169-007 Porto, Portugal}\\
{\em email:} pvsilva@fc.up.pt\\
$ $\\
{\bf Mihalis Sykiotis}\\ $ $\\ {\em Department of Mathematics,
  National and Kapodistrian University of Athens,}\\
{\em Panepistimioupolis, GR-157 84, Athens, Greece}\\
{\em e-mail:} msykiot@math.uoa.gr}

\date{\today}

\begin{document}
\maketitle

\begin{center}\small

\bigskip

\end{center}

\abstract{We discuss in the context of finite extensions two classical
  theorems of Takahasi and Howson on subgroups of free groups. We
  provide bounds for the rank of the intersection of subgroups within
  classes of groups such as virtually free groups, virtually nilpotent
  groups or fundamental groups of finite graphs of groups with
virtually polycyclic vertex groups and finite edge groups. As an
application of our generalization of Takahasi's Theorem, we provide an
uniform bound for the rank of the periodic subgroup of any
endomorphism of the fundamental group of a given finite graph of groups with
finitely generated virtually nilpotent vertex groups and finite edge
groups.}

\section{Introduction}

Some famous theorems on subgroups of free groups involve finiteness
conditions. Part of them admit generalizations to further classes of
groups, and constructions such as free products, finite extensions or
graphs of groups have been involved in most of them.

For instance, Howson's Theorem states that the intersection
of two finitely generated subgroups $H,K$ of a free group is also finitely
generated. In his seminal paper \cite{How}, Howson also provided an
upper bound on the
rank of $H \cap K$ with respect to the ranks of $H$ and $K$, namely
(for $H$ and $K$ nontrivial):
$$\rk(H \cap K) \leq 2\rk(H)\rk(K) - \rk(H) - \rk(K) +1.$$
Later on, Hanna Neumann improved this upper bound to
$$\rk(H \cap K) \leq 2(\rk(H) -1)(\rk(K)-1) +1$$
and conjectured that the factor 2 could be removed, the famous {\em
  Hanna Neumann Conjecture}. The Conjecture was finally proved in 2011
by Friedman and Mineyev (independently):

\bt
\label{mine}
{\rm \cite{Fri,Min}}
Let $F$ be a free group and let $K_1,K_2 \leq F$ be finitely generated and
nontrivial. Then
$${\rm rk}(K_1\cap K_2) \leq ({\rm rk}(K_1) -1)({\rm rk}(K_2) -1) +1.$$
\et

Howson's Theorem led to the concept of Howson group: a group $G$ is
a {\em Howson group} if the intersection of finitely generated
subgroups of $G$ is still finitely generated. Kapovich has shown
that many hyperbolic groups fail this property \cite{Kap}, but it is
easy to show that Howson groups are closed under finite extension
and so in particular virtually free groups are Howson groups. More
generally, the class of Howson groups is closed under graphs of
groups, where the edge groups are finite (see \cite[Theorem 2.13
(1)]{Syk3} for a proof). But can we get some rank formula as in the
case of free groups? A recent paper of Zakharov \cite{Zak} provides
an upper bound for the rank of the intersection of two {\em free}
finitely generated subgroups of a virtually free group. In the case
of free products, upper bounds for the Kurosh rank of the
intersection of subgroups have been obtained by various authors. See
for instance \cite{AMS}, and the references therein, where Theorem
\ref{mine} is extended to free products of right-orderable groups.

We introduce in Section \ref{sh} the concept of {\em strongly Howson
  group}, when an uniform bound for the rank of $H \cap K$ can be
obtained from bounds on the ranks of $H$ and $K$. We show that the
class of strongly Howson groups is closed under finite extensions
and compute bounds using an improved version of Schreier's Lemma,
which can be obtained with the help of Stallings automata. These
bounds are then applied to several particular cases such as
virtually free, virtually polycyclic, virtually nilpotent, and more
generally fundamental groups of finite graphs of groups with
virtually polycyclic vertex groups and finite edge groups.

Another famous result, known as Takahasi's Theorem, states the following:

\bt
\label{taka}
{\rm \cite{Tak}}
Let $F$ be a free group and let
$K_1 \leq K_2 \leq \ldots$ be an ascending chain of finitely generated
subgroups of $F$. If the rank of the subgroups in the chain is
bounded, then the chain is stationary.
\et

Bogopolski and Bux proved recently an analogue of Takahasi's Theorem
for fundamental groups of closed compact surfaces \cite[Proposition
2.2]{BB}. We say that a group $G$ is a {\em Takahasi group} if every
ascending chain $H_{1}\leq H_{2}\leq \cdots$  of subgroups each of
rank $\leq M$ in $G$, is stationary. We prove, in Section \ref{st},
that the class of Takahasi groups is closed under finite extensions
and finite graphs of groups with virtually polycyclic vertex groups
and finite edge groups.

We provide an application of the generalized Takahasi's Theorem in
Section \ref{spp}. Finally, using previous work of the third author
\cite{Syk}, we show that the periodic subgroup is finitely generated
for every endomorphism of the fundamental group of a finite graph of
groups with finitely generated virtually nilpotent vertex groups and
finite edge groups. As a consequence, we can bound the periods for
each particular endomorphism of such a group.

\section{Preliminaries}

We collect in this section some standard group-theoretic concepts and
results. The reader is referred to \cite{Hal,KM,LS} for details.

Given a group $G$ and $X \subseteq G$, we denote by $\langle X
\rangle$ the subgroup of $G$ generated by $X$. If $G$ is finitely
generated, the {\em rank} of $G$ is defined as
$$\rk(G) = \min\{ |X| : G = \langle X \rangle \}.$$

We denote by $F_A$ the free group on an alphabet $A$. A free group of
rank $n$ is generically denoted by $F_n$. The standard way of
describing finitely generated subgroups of a free group is by means
of Stallings automata, a
construction designed by Stallings under a different formalism
\cite{Sta2}.

To simplify things, we define an {\em automaton} to be a
structure of the form $\A = (A,Q,q_0,T,E)$ where:
\begin{itemize}
\item
$A$ is a finite alphabet;
\item
$Q$ is a set (vertices);
\item
$q_0 \in Q$ (initial vertex);
\item
$T \subseteq Q$ (terminal vertices);
\item
$E \subseteq Q \times A \times Q$ (edges).
\end{itemize}
The automaton is finite if $Q$ is finite.

A {\em finite nontrivial path} in $\A$ is a sequence
$$p_0 \mapright{a_1} p_1 \mapright{a_2} \ldots \mapright{a_n} p_n$$
with $(p_{i-1},a_i,p_i) \in E$ for $i = 1,\ldots,n$. Its {\em label}
is the word $a_1\ldots a_n \in A^*$. It is said to be a {\em
  successful} path if $p_0 = q_0$ and $p_n \in T$. We consider also
the {\em trivial path} $p \mapright{1} p$ for $p \in Q$. It is
successful if $p = q_0 \in T$.

The {\em language} $L(\A)$ {\em
  recognized by} $\A$
is the set of all labels of successful paths in $\A$. For
details on automata, the reader is referred to \cite{Ber,Sak}.

Let $H \leq
F_A$ be finitely generated. Taking a finite set of generators
$h_1,\ldots,h_n$ of $H$ in
reduced form, we start with the so-called flower automaton $\F(H)$ (on
the alphabet $A \cup A\inv$), where
{\em petals} (of variable length) labelled by the words $h_i$ are
glued to a basepoint $q_0$ (both initial and
terminal):
$$
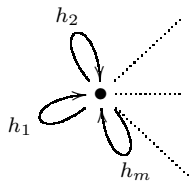
\xymatrix{
& \ \\
\bullet
\ar@(dl,l)^{h_1}
\ar@(ul,u)^{h_2}
\ar@{.}[ur] \ar@{.}[r] \ar@{.}[dr] \ar@(dr,d)^{h_m} & \ \\
& \ }$$
We include also an edge of the form $q \mapright{a\inv} p$ for every
edge of the form $p \mapright{a} q$.
Then we proceed by successively folding pairs of edges of the form $q
\mapleft{a} p \mapright{a} r$ $(a \in A \cup A\inv)$.
The final automaton $\S(H)$ does not depend on the folding sequence
nor even on the original finite generating set, and is known as the
{\em Stallings automaton} of $H$. For details and
applications of Stallings automata, see \cite{BS, KM, MVW}.

One of the classical applications of Stallings automata provides a
solution for the generalized word problem of $F_A$ (see
\cite[Proposition 2.5]{BS}): given $u \in F_A$ in reduced form, we
have
\beq
\label{gwp}
u \in H \hspace{.3cm} \iff  \hspace{.3cm} u \in L(\S(H)).
\eeq
Another famous application (see \cite[Proposition 2.6]{BS}) is
the rank formula
$$\rk(H) = e - v +1,$$
where $v$ denotes the number of vertices of $\S(H)$ and $e$ denotes
the number of positive edges of $\S(H)$ (i.e. edges labelled by
letters of $A$). In the particular case where $H$ is of finite index
in $F_A$, we get $v = [F_A:H]$ and $e = [F_A:H]|A|$, hence (see
\cite[Proposition 2.8]{BS}) \beq \label{rsa} \rk(H) = [F_A:H](|A|-1)
+1. \eeq

Given a class $\C$ of groups, we say that a group $G$ is:
\bi
\item {\em virtually} $\C$ if $G$ has a finite index subgroup in $\C$;
\item $\C$-{\em by-finite} if $G$ has a finite index normal subgroup in $\C$.
\ei

If the class $\C$ is closed under isomorphism and taking subgroups, then the
two concepts coincide. That is the case for free, nilpotent, polycyclic and
strongly polycyclic groups.

If $F$ is a finite index subgroup of
$G$, we also say that $G$ is a {\em finite extension} of $F$.
If $[G:F] = m$, we
may decompose $G$ as a disjoint union of right cosets
\beq
\label{stdec}
G = Fb_1 \cup \ldots \cup Fb_m
\eeq
with $b_1 = 1$. We shall refer to (\ref{stdec}) as a {\em standard
  decomposition} of $G$ with respect to $F$.

The next simple result is essential to handle subgroups of finite
extensions:

\bp
\label{sdsub}
Let $G$ be a finite extension of a group $F$ with standard decomposition
(\ref{stdec}). Let $H \leq G$ and write $K = H \cap F$. Then there
exist $I \subseteq \{ 2, \ldots,m\}$ and $x_i \in F$ $(i \in I)$ such
that
\beq
\label{sdsub1}
H = K \cup (\bigcup_{i \in I} Kx_ib_i).
\eeq
\ep

\proof
Let
$$I = \{ i \in \{ 2, \ldots,m\} \mid H \cap Fb_i \neq \emptyset\}.$$
Since $b_1 = 1$, we may write
$$H = K \cup (\bigcup_{i \in I} K_ib_i)$$
for some nonempty $K_i \subseteq F$ $(i \in I)$. For each $i \in I$,
fix $x_i \in K_i$. It remains to be proved that $K_i = Kx_i$.

Clearly, $Kx_ib_i \subseteq HK_ib_i \subseteq H^2 = H$, hence $Kx_i
\subseteq K_i$. Conversely, let $y \in K_i$. Then $yx_i\inv =
(yb_i)(x_ib)\inv \in HH\inv = H$. Since also $yx_i\inv \in K_iK_i\inv
\subseteq FF\inv = F$, we get $yx_i\inv \in K$ and so $y \in
Kx_i$. Thus $K_i = Kx_i$ as required.
\qed


Now we recall the definitions of several other classes of groups which
play a part in this paper.

A group $G$ is {\em residually finite} if the intersection of all
normal subgroups of finite index is equal to ${1}$. Since any
subgroup of finite index in a group $G$ contains a normal subgroup
of finite index in $G$, it follows that a finite extension of a
residually finite group is residually finite.


Let $G$ be a group. Given $H,K \leq G$, write
$$[H,K] = \langle hkh\inv k\inv : h \in H\; k \in K \rangle \leq G.$$
The {\em lower central series} of $G$ is the
sequence
$$G = G_0 \unrhd G_1 \unrhd G_2 \unrhd \ldots,$$
where $G_{n} = [G,G_{n-1}]$ for every $n \geq 1$. The group $G$ is
{\em nilpotent} if $G_n = \{ 1 \}$ for some $n \geq 1$. The minimum
such $n$ is the {\em nilpotency class} of $G$. Clearly, an abelian
group is nilpotent of nilpotency class $\leq 1$. A subgroup of a
nilpotent group of class $n$ is nilpotent of class $\leq n$.



A group $G$ is called {\em polycyclic} if it admits a subnormal series
\beq
\label{ycy}
G = G_0 \rhd G_1 \rhd \ldots \rhd G_n = \{ 1 \}
\eeq
such that $G_{i-1}/G_i$ is cyclic for $i = 1,\ldots,n$. The minimum
such $n$ is the {\em polycyclic rank} of $G$ and is denoted by
$\prk(G)$. The {\em
  Hirsch number} $h(G)$ of $G$ is defined as the number of infinite
factors $G_{i-1}/G_i$, which is independent from the subnormal
series. In particular, $h(G) \leq \prk(G)$. If $G_{i-1}/G_i$ is
infinite cyclic for every $i$, we say that $G$ is {\em strongly
polycyclic}. Every polycyclic group has a normal strongly polycyclic
group of finite index. By \cite{Hir}, polycyclic groups are
residually finite.

The class of (strongly) polycyclic groups is closed under taking
subgroups. Moreover, a simple induction on $\prk(G)$ shows that if $G$
is polycyclic then
\beq
\label{shpol2}
\rk(H) \leq \prk(G) \hspace{.5cm}\mbox{for every }H \leq G.
\eeq

A nilpotent
group is polycyclic if and only if it is finitely
generated. For details on polycyclic groups, see \cite{Weh}.

Finally, we recall the concept of graph of groups, central in Bass-Serre
theory \cite{Ser}.

Following Serre, a {\em graph} is a
structure of the form
$\Gamma = (V,E,\alpha, \, \bar{}\, )$, where:
\bi
\item
$V$ is a nonempty set (vertices);
\item
$E$ is a set (edges);
\item
$\alpha:E \to V$ is a mapping;
\item
$\bar{}:E \to E$ is an involution without fixed points.
\ei
Concepts such as cycle, connectedness, tree or
subgraph are defined the obvious way.
If $\Gamma$ is connected and $T \subseteq E$ defines a subtree of
$\Gamma$ connecting all the vertices, we say that $T$ is a {\em
  spanning tree} of $\Gamma$.

A (finite) {\em graph of groups} over a (finite) connected graph
$\Gamma$ is a structure of the form
\beq
\label{gg}
\G = ((G_v)_{v \in V}, (G_e)_{e \in E}, (\alpha_e)_{e \in E}),
\eeq
where:
\bi
\item
the $G_v$ are groups for all $v \in V$ (vertex groups);
\item
the $G_e$ are groups for all $e \in E$ (edge groups) satisfying
$G_{\bar{e}} = G_e$;
\item
the $\alpha_e:G_e \to G_{e\alpha}$ are monomorphisms for all $e \in E$
(boundary monomorphisms).
\ei

The {\em fundamental group} $\pi_1(\G,T)$ of the graph of groups
(\ref{gg}) with respect to a spanning tree $T$ of $\Gamma$ is the
quotient of the free product
$$(\ast_{v \in V} G_v) \ast F_E$$
by the normal subgroup generated by the following elements:
\bi
\item
$e\bar{e}$ $(e \in E)$;
\item
$t \in T$;
\item
$e\inv(g\alpha_e)e(g\alpha_{\bar{e}})\inv$ $(e \in E,\; g \in G_e)$.
\ei
The vertex groups are naturally embedded into $\pi_1(\G,T)$, which is
independent of the chosen spanning tree $T$, up to isomorphism.

If the edge groups $G_e$ are all trivial, then we get a free
product
\beq
\label{egt}
\pi_1(\G,T) = (\ast_{v \in V} G_v) \ast F_{A},
\eeq
where $E \setminus T = A \cup \bar{A}$ and $A \cap A\inv = \emptyset$.
HNN extensions and
amalgamated free products constitute important particular cases of
this construction, by taking graphs with two edges, respectively of the form
$$
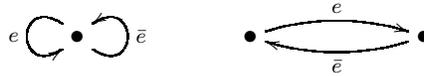
\xymatrix{
\bullet \ar@(ul,dl)_e \ar@(dr,ur)_{\bar{e}} && \bullet \ar@/^/[rr]^e
&& \bullet \ar@/^/[ll]^{\bar{e}}
}$$
Moreover, whenever $\Gamma$ is finite, the fundamental group
$\pi_1(\G,T)$ can be built from the vertex groups
using a finite number of HNN extensions and
amalgamated free products, where the associated/amal\-ga\-mated subgroups
are of the form $G_e\alpha_e$.

The nature of $\pi_1(\G,T)$ is conditioned by the nature of the
vertex and edge groups. This is illustrated by the following
well-known theorem of Karrass, Pietrowski and Solitar \cite{KPS2} (see
also \cite[Theorem 7.3]{SW}): a finitely generated group is virtually free if
and only if it is the fundamental group of a finite graph of finite
groups.

\section{Howson's Theorem}
\label{sh}

We say that a group $G$ is {\em strongly Howson} if
\beq
\label{strhow}
\sup \{ \rk(H_1\cap H_2) \mid H_1,H_2 \leq G,\; \rk(H_1) \leq n_1,\;
\rk(H_2) \leq n_2\} < \infty
\eeq
for all $n_1,n_2 \in \N$. In this case, we can define a function $\xi_G: \N
\times \N \to \N$ by letting (\ref{strhow}) be $(n_1,n_2)\xi_G$. Since
every subgroup of a cyclic group is cyclic, we only care about the
nontrivial cases $n_1,n_2 \geq 2$.

Clearly,
if $G$ is strongly Howson and $H \leq G$, then $H$ is strongly Howson
and $(n_1,n_2)\xi_H \leq (n_1,n_2)\xi_G$ for all $n_1,n_2 \in \N$.

Trivially, every strongly Howson group is a Howson group. We ignore if the
converse is true.

Schreier's Lemma \cite{Sch} (see also \cite{Hal}) states that the
inequality
$${\rm rk}(H) \leq [G:H]{\rm rk}(G)$$
holds whenever $H$ is a finite index subgroup of a finitely
generated group $G$. The following improved version is well-known,
but we give a proof for completeness, using Stallings automata:

\bp
\label{newrankfi}
Let $H$ be a finite index subgroup of a finitely generated group
$G$. Then ${\rm rk}(H) \leq [G:H]({\rm rk}(G) - 1) +1$.
\ep

\proof
Let $m = [G:H]$ and $n = \rk(G)$. Then there exists an epimorphism
$\p:F_n \to G$. It is straightforward that $[F_n:H\p\inv] = [G:H] =
m$, hence it follows from (\ref{rsa}) that
$$\rk(H\p\inv) = [F_n:H\p\inv](n-1) +1 = m(n-1) +1.$$
Since $H = (H\p\inv)\p$ yields $\rk(H) \leq \rk(H\p\inv)$, we are
done.
\qed

Note that this bound is tight in view of (\ref{rsa}).

Now we can prove the following result:

\bt
\label{boho}
Let $G$ be a finite extension of a strongly Howson group $F$ and let
$m = [G:F]$. Then $G$ is strongly Howson and
$$(n_1,n_2)\xi_G \leq (m(n_1-1)+1,m(n_2-1)+1)\xi_F +m-1$$
for all $n_1,n_2 \geq 1$.
\et

\proof
Let $H_1,H_2 \leq G$ with $\rk(H_j) \leq n_j$ for $j = 1,2$.
We must show that
\beq
\label{boho5}
\rk(H_1 \cap H_2) \leq (m(n_1-1)+1,m(n_2-1)+1)\xi_F +m-1.
\eeq
By Proposition \ref{sdsub}, for $j = 1,2$ we may write
$$H_j = K_j \cup (\bigcup_{i \in I_j} K_jx_i^{(j)}b_i)$$
with $K_j = H_j \cap F$, $I_j \subseteq \{ 2, \ldots,m\}$ and
$x_i^{(j)} \in F$ $(i \in I_j)$.

For all $h,h' \in H_j$,
$$Fh = Fh' \Rw h'h\inv \in F \Rw h'h\inv \in H_j \cap F = K_j \Rw K_jh = K_jh',$$
hence
$$[H_j:K_j] \leq [G:F] = m$$
and Proposition \ref{newrankfi} yields
\beq
\label{boho2}
\rk(K_j) \leq m(\rk(H_j)-1)+1 \leq m(n_j-1)+1.
\eeq
On the other hand, writing $K = K_1 \cap K_2 = H_1 \cap H_2 \cap F$,
it follows from Proposition \ref{sdsub} that
$$H_1\cap H_2 = K \cup (\bigcup_{i \in I} Ky_ib_i)$$
for some $I \subseteq \{ 2, \ldots,m\}$ and
$y_i \in F$ $(i \in I)$. Since
$$H_1\cap H_2 = \langle K \cup \{ y_ib_i \mid i \in I \} \rangle,$$
we get
\beq
\label{boho3}
\rk(H_1 \cap H_2) \leq \rk(K) + |I| \leq \rk(K) + m-1.
\eeq
In view of (\ref{boho2}), we get
$$\rk(K) \leq (m(n_1-1)+1,m(n_2-1)+1)\xi_F$$
and so
(\ref{boho3}) yields
$$\rk(H_1 \cap H_2) \leq (m(n_1-1)+1,m(n_2-1)+1)\xi_F + m-1.$$
Therefore (\ref{boho5}) holds as required.
\qed

We apply Theorem \ref{boho} to some classes of groups, starting with
the straightforward virtually free case:

\bc
\label{showvf}
Let $G$ be a virtually free group with a free subgroup $F$ of index
$m$. Then $G$ is strongly Howson and
$$(n_1,n_2)\xi_G \leq m^2(n_1-1)(n_2-1) +m$$
for all $n_1,n_2 \geq 1$.
\ec

\proof
By Theorem \ref{mine}, we have
$$(k_1,k_2)\xi_F \leq (k_1-1)(k_2-1)+1$$
for all $k_1,k_2 \geq 1$. By Theorem \ref{boho}, we get
$$\begin{array}{lll}
(n_1,n_2)\xi_G&\leq&(m(n_1-1)+1,m(n_2-1)+1)\xi_F +m-1
\leq (m(n_1-1))(m(n_2-1)) +m\\
&=&m^2(n_1-1)(n_2-1) +m
\end{array}$$
In particular, $G$ is strongly Howson.
\qed

In a recent paper, Zakharov proved the following theorem:

\bt
\label{zak}
{\rm \cite[Theorem 2]{Zak}}
Let $G$ be a virtually free group and let
$H_1,H_2 \leq G$ be finitely generated, free and nontrivial. Then
$${\rm rk}(H_1\cap H_2) \leq 6n({\rm rk}(H_1) -1)({\rm rk}(H_2) -1)
+ 1,$$
where $n$ is the maximum of orders $|P \cap (H_1H_2)|$ over all finite
subgroups $P$ of $G$. As a consequence,
$${\rm rk}(H_1\cap H_2) \leq 6m({\rm rk}(H_1) -1)({\rm rk}(H_2) -1)
+ 1$$
if $G$ has a free subgroup of index $m$.
\et

How does the upper bound arising from Corolllary \ref{showvf} compare
with the upper bounds in
Theorem \ref{zak}? In general, for arbitrary free subgroups, the
bounds in Theorem \ref{zak} are
smaller since they are linear on $m$ and ours are quadratic. However,
we claim that our bound is actually smaller than the second bound in
Theorem \ref{zak} if $m \leq 5$ and
$H_1,H_2$ are noncyclic (if $H_1$ or $H_2$ is cyclic, so is $H_1 \cap
H_2$ and we have a trivial case anyway). Indeed, the product $p =
(\rk(H_1)-1)(\rk(H_2)-1)$ is then positive and so
$$\begin{array}{ll}
&m^2(\rk(H_1)-1)(\rk(H_2)-1) + m < 6m(\rk(H_1)-1)(\rk(H_2)-1) + 1\\
\iff&m^2p + m \leq 6mp \iff mp + 1 \leq 6p \iff mp < 6p \iff m < 6.
\end{array}$$

The following example shows that our bound may also beat the first
bound provided by Theorem \ref{zak}:

\be
Let $A = \{ a,b,c \}$ and let $C_2$
be a cyclic group of order 2. Let $G = F_A \times C_2$ and
$$H_1 = \langle (a,1),(bc,1)\rangle, \quad H_2 = \langle (ab,1),(c,0)\rangle.$$
Then:
\bi
\item[(i)] $H_1$ and $H_2$ are free subgroups of rank 2 of the virtually free
  group $G$;
\item[(ii)] Theorem \ref{zak} provides the upper bound ${\rm
    rk}(H_1\cap H_2) \leq 13$;
\item[(iii)] Theorem \ref{boho} provides the upper bound ${\rm
    rk}(H_1\cap H_2) \leq 6$;
\item[(iii)] actually, ${\rm
    rk}(H_1\cap H_2) = 1$.
\ei
\ee

Indeed, $F_A \times \{ 0 \}$ is a free subgroup of index 2 of $G$, hence
$G$ is virtually free.

It is easy to see that projecting $H_1$ into its first
component we get a free group with basis $\{ a,bc \}$, and we can deduce
from that fact that $H_1$ is itself free of rank 2. Similarly, $H_2$ is
free of rank 2.

Let $P = \{ 1 \} \times C_2 \leq G$. It is easy to check that $|P \cap
H_1H_2| = 2$, e.g.
$$(1,1) = (a,1)(bc,1)(c,0)\inv (ab,1)\inv \in H_1H_2$$
and so we get the upper bound
$\rk(H_1\cap H_2) \leq 13$ from Theorem \ref{zak}. On the other hand,
it is immediate that
Corollary \ref{showvf} yields the upper bound
$\rk(H_1\cap H_2) \leq 6$.

Finally, with the help of the standard algorithm to compute a basis for
the intersection in free groups \cite[Proposition 9.4]{KM}, it is easy
to check that
$$\langle a,bc\rangle \cap \langle ab, c\rangle = \langle
abc\rangle.$$
It follows easily that
$$H_1 \cap H_2 = \langle ((abc)^2,0)\rangle$$ and so $\rk(H_1\cap H_2) =
1$.

\medskip

We present further applications of Theorem \ref{boho}:

\bc
\label{shpol}
Let $G$ be a virtually polycyclic group. Then $G$ is strongly Howson
and $\xi_G$ is a bounded function.
\ec

\proof
Let $P$ be a polycyclic subgroup of $G$ of index $m$. Let $n = \prk(G)$.
%
%
By (\ref{shpol2}), we have
$$(n_1,n_2)\xi_P \leq n$$
for all $n_1,n_2 \in \N$. By Theorem
\ref{boho}, we get
$(n_1,n_2)\xi_G \leq n+m-1$ for all $n_1,n_2 \geq 1$.
Thus $\xi_G$ is bounded and $G$ is strongly Howson.
\qed

The general virtually nilpotent case is a bit harder. Note that a non
finitely generated nilpotent group is not polycyclic.

\bt
\label{shnil}
Let $G$ be a virtually nilpotent group. Then $G$ is strongly Howson and
$$(n_1,n_2)\xi_G \leq
  \frac{(m(p-1)+1)^{n+1}-m(p-1)-1}{m(p-1)} +m-1$$
for all $n_1,n_2 \geq 2$
  and $p = {\rm min}\{ n_1,n_2\}$.
\et

\proof
Suppose that $N$ is a nilpotent group of class $n$ and rank
$k \geq 2$. We claim that
\beq
\label{shnil1}
\rk(H) \leq \frac{k^{n+1}-k}{k-1}
\eeq
for every $H \leq N$.

Let \beq \label{shnil4} N = N_0 \rhd N_1 \rhd \ldots \rhd N_n = \{ 1
\} \eeq be the lower central series of $N$. By \cite[Corollary
10.3]{Hal}, we have \beq \label{shnil2} \rk(N_{i-1}/N_i) \leq k^i
\eeq for $i = 1,\ldots,n$. Since $[N_{i-1},N_{i-1}] \subseteq
[N,N_{i-1}] = N_i$, the quotient $N_{i-1}/N_i$ is abelian. By
(\ref{shnil2}), there exist $x_1, \ldots, x_{k^i} \in N_{i-1}$ such
that
$$N_{i-1}/N_i = \langle x_1N_i, \ldots x_{k^i}N_i \rangle.$$
Let $\pi_i:N_{i-1} \to N_{i-1}/N_i$ be the canonical
projection. For $j = 0, \ldots, k^i$, let
$$N_{i,j} = \langle x_1N_i, \ldots x_{j}N_i \rangle\pi_i\inv.$$
Since $\langle x_1N_i, \ldots x_{j-1}N_i \rangle \unlhd \langle x_1N_i,
\ldots x_{j}N_i \rangle$ due to $N_{i-1}/N_i$ being abelian, we get
$N_{i,j-1} \unlhd N_{i,j}$ and so we have a chain
\beq
\label{shnil3}
N_{i-1} = N_{i,k^i} \unrhd \ldots \unrhd N_{i,1} \unrhd N_{i,0} =
N_i.
\eeq
Moreover,
$$N_{i,j}/N_{i,j-1} = \langle x_1N_i, \ldots x_{j}N_i
\rangle\pi_i\inv\, / \, \langle x_1N_i, \ldots x_{j-1}N_i \rangle\pi_i\inv
\cong \langle x_1N_i, \ldots x_{j}N_i
\rangle / \langle x_1N_i, \ldots x_{j-1}N_i
\rangle$$
and is therefore cyclic since $N_{i-1}/N_i$ is abelian.

Inserting the chains (\ref{shnil3}) into (\ref{shnil4}), we obtain
a subnormal series for $N$ with length
$$k + k^2 + \ldots + k^n = \frac{k^{n+1}-k}{k-1}$$
and cyclic quotients. In particular, $N$ is polycyclic. Now
(\ref{shnil1}) follows from (\ref{shpol2}).

Assume now that $N$ is a nilpotent subgroup of $G$ of class $n$ and
index $m$. Let $n_1,n_2 \geq 2$ and suppose that $H_1,H_2 \leq N$ are
such that $\rk(H_j) \leq n_j$ for $j = 1,2$. Since each $H_j$ is also
nilpotent of class $\leq n$ and $H_1 \cap H_2 \leq H_j$,
(\ref{shnil1}) yields
$$\rk(H_1 \cap H_2) \leq \frac{n_j^{n+1}-n_j}{n_j-1}$$
and so, writing $p = \min \{ n_1, n_2 \}$, we get
$$(n_1,n_2)\xi_N \leq p + p^2 + \ldots + p^n = \frac{p^{n+1}-p}{p-1}.$$

In particular, $N$ is strongly Howson and we may apply Theorem
\ref{boho} to get
$$(n_1,n_2)\xi_G \leq
(m(n_1-1)+1,m(n_2-1)+1)\xi_N +m-1.$$
Since $\min\{m(n_1-1)+1,m(n_2-1)+1 \} = m(p-1)+1$, we get
$$(n_1,n_2)\xi_G \leq \frac{(m(p-1)+1)^{n+1}-m(p-1)-1}{m(p-1)} +m-1.$$
\qed

Our last application involves graphs of groups, but first we deal with
the following particular case:

\bt
\label{shspg}
Let $G = S_1 \ast \ldots \ast S_t$ be a free product of strongly
polycyclic groups and let $M = {\rm max}\{ h(S_1), \ldots, h(S_t)
\}$. Then $G$ is strongly Howson and
$$(n_1,n_2)\xi_G \leq M(n_1-1)(n_2-1)+M$$
for all $n_1,n_2 \geq 1$.
\et

\proof
In view of (\ref{shpol2}), we have
$\rk(L) \leq M$
for all $i \in \{ 1, \ldots,t\}$ and $L \leq S_i$.

By the Kurosh subgroup theorem,
every subgroup $H \leq G$ is isomorphic to a free product of the form
$$(\ast_{j \in J} L_j) \ast F_{A},$$
where each $L_j$ is the intersection of $H$ with some conjugate of some
$S_i$. The {\em Kurosh rank} of $H$ is defined by
$$\krk(H) = |J|+|A|.$$
It follows from Grushko Theorem on the additivity of ranks in free
products \cite{Gru} that \beq \label{utah} \krk(H) \leq \rk(H). \eeq
In general finite Kurosh rank does not imply finite rank. But in the
present case, since $\rk(L_j) \leq M$ for every $j \in J$, we have
\beq \label{bsvp1} \rk(H) \leq \sum_{j \in J} \rk(L_j) + |A| \leq
M|J| + |A| \leq M\krk(H). \eeq

Let $H_j \leq G$ satisfy $\rk(H_j) \leq n_j$ for $j = 1,2$. We may
assume that $H_1 \cap H_2$ is nontrivial. Now (\ref{utah}) yields
$\krk(H_j) \leq \rk(H_j) \leq n_j$. Since strongly polycyclic groups
are right-orderable \cite{Rhe},
it follows from \cite[Theorem A]{AMS}
that
$$\krk(H_1 \cap H_2) \leq (\krk(H_1) -1)(\krk(H_2) -1)+1.$$
Hence (\ref{utah}) and (\ref{bsvp1}) yield
$$\begin{array}{lll}
\rk(H_1 \cap H_2)&\leq&M\krk(H_1 \cap H_2) \leq M(\krk(H_1)
-1)(\krk(H_2) -1)+M\\
&\leq&M(\rk(H_1)
-1)(\rk(H_2) -1)+M \leq M(n_1-1)(n_2-1)+M.
\end{array}$$
\qed

Now we prove the following lemma:

\bl
\label{find}
Let $G$ be the fundamental group of a finite graph of groups $\G$ with
finite edge groups.
\bi
\item[(i)] If $\G$ has virtually polycyclic vertex groups, then $G$
has a finite index normal subgroup which is a finitary free product of strongly
polycyclic groups.
\item[(ii)] If $\G$ has finitely generated virtually nilpotent vertex
  groups, then $G$ has a finite index normal subgroup which is a
  finitary free product of finitely generated nilpotent groups.
\ei
\el

\proof
(i) Let
$$\G = ((G_v)_{v \in V}, (G_e)_{e \in E}, (\alpha_e)_{e \in E})$$
be such a graph of groups, built over the finite connected graph
$$\Gamma = (V,E,\alpha, \, \bar{}\, ).$$
Fix a spanning tree $T$ of $\Gamma$ and let $G = \pi_1(\G,T)$. Since
polycyclic groups are residually finite, it follows that each vertex
group $G_v$ is residually finite. Now the class of residually finite
groups is closed under amalgamated free products with finite
amalgamated subgroups and under HNN extensions with finite
associated subgroups \cite{Bau,Coh}. Since $\Gamma$ is a finite
graph, we may use the decomposition of $G$ in terms of HNN
extensions and amalgamated products over the finite edge groups to
deduce that $G$ is itself residually finite.

Let
$$X = (\bigcup_{e \in E} G_e\alpha_e)\setminus \{ 1 \} \subseteq G
\setminus \{ 1 \}$$
consist of the image of the edge groups in
$G$ through the boundary monomorphisms, with the identity
removed. Since both the graph and
the edge groups are finite, so is $X$. Let $x \in X$. Since $G$ is
residually finite, there exists some $N_x \lhd G$ of finite index such
that $x \notin N_x$. Let
$$N = \bigcap_{x \in X} N_x.$$
Since $X$ is finite, $N$ is still a normal subgroup of $G$ of finite
index.

By \cite[Corollary 2]{KPS}, since $G$ is the fundamental group of a
finite graph of groups
with finite edge groups, every finite index $H \leq G$ is itself
the fundamental group of a
finite graph of groups $\G_H$ where:
\bi
\item
the vertex groups are conjugates of subgroups of the form $H \cap
yG_vy\inv$ $(v \in V,\; y \in G)$;
\item
the edge groups are conjugates of subgroups of the form $H \cap
y(G_e\alpha_e)y\inv$ $(e \in E,\; y \in G)$.
\ei

We consider now the case $H = N$. Since $N \cap X = \emptyset$ by
construction, we have $N \cap G_e\alpha_e = \{ 1
\}$ for every $e \in E$. Since $N$ is normal, we get
$$N \cap y(G_e\alpha_e)y\inv = y(y\inv Ny \cap G_e\alpha_e)y\inv =
y(N \cap G_e\alpha_e)y\inv = 1.$$ Thus $\G_N$ has trivial edge
groups.

On the other hand, if $G'$ has a polycyclic subgroup $F'$ of index
$m$ and $H' \leq G'$, it follows from Proposition \ref{sdsub} that
$[H':H' \cap F'] \leq m$. Since $H' \cap F'$ must be itself
polycyclic, $H'$ is virtually polycyclic as well.

Thus each
group $H \cap
yG_vy\inv$ is virtually polycyclic and so $\G_N$ has virtually
polycyclic vertex groups.

But then $N$, being the fundamental group of $\G_N$, is a free
product of finitely many virtually polycyclic groups and a free
group of finite rank. Since a free group of finite rank is the free
product of finitely many cyclic groups (hence polycyclic), it
follows that $N$ is indeed a free product of finitely many virtually
polycyclic groups, say $N = K_1 \ast \ldots \ast K_t$, with the
$K_i$ nontrivial.

Since $K_i$ is indeed virtually strongly polycyclic for each $i$, it
contains a strongly polycyclic subgroup $P_i$ of finite index for $i
= 1,\ldots,t$. Since a subgroup of $P_i$ must be still strongly
polycyclic, we may assume that $P_i \unlhd K_i$. Let
$$\p:N \to K_1/P_1 \times \ldots \times K_t/P_t$$
be the canonical epimorphism. Then $\ker(\p)$ is a finite index
subgroup of $N$. Since $[G:N] < \infty$, we have $[G:\ker(\p)] < \infty$ as
well. Let
$$Q = \bigcap_{g \in G} g(\ker(\p))g\inv.$$
Since $[G:\ker(\p)] < \infty$, $Q$ is a finite index normal subgroup
of $G$.

Since $Q \subseteq N = K_1 \ast \ldots
\ast K_t$, it follows from the Kurosh subgroup theorem \cite{Kur} that
$Q$ is isomorphic to a free product of the form
$$(\ast_{j \in J} L_j) \ast F_A,$$
where each $L_j$ is the intersection of $Q$ with some conjugate of some
$K_i$. Now
$L_j = Q \cap y_jK_iy_j\inv$ implies
$$L_j = y_j(y_j\inv Q y_j \cap K_i)y_j\inv = y_j(Q \cap K_i)y_j\inv
\subseteq y_j(\ker(\p) \cap K_i)y_j\inv \subseteq y_jP_iy_j\inv$$
and so $L_j$, being a subgroup of a strongly polycyclic group, is also
strongly polycyclic. Since $F_A$ is a free product of cyclic groups,
it follows that $Q$ is a finite index normal subgroup of $G$ which is a free
product of strongly polycyclic groups.

(ii) Each vertex group is a finite extension of a finitely generated
nilpotent group, therefore the vertex groups are virtually
polycyclic. Thus we only need to perform minimal adaptations to the
proof of (i) which we proceed to enhance: \bi
\item
Since the class of nilpotent groups is
closed under taking subgroups, the same happens with the class of
finitely generated nilpotent groups (since they are polycyclic and in
view of (\ref{shpol2})) and therefore with the
class of finitely generated virtually nilpotent groups (in view of
Proposition \ref{sdsub}). Thus $\G_N$ has finitely generated virtually
nilpotent vertex groups.
\item
$N$ is the free product of finitely many finitely generated virtually
nilpotent groups and a free group of finite rank. Since $\Z$ is
nilpotent, then $N$ is the free product of finitely many finitely
generated virtually nilpotent groups.
\item
We choose the $P_i$ to be finitely generated nilpotent. The free
factors of $Q$ are then finitely generated nilpotent groups.
\ei
\qed

We can finally prove the following:

\bt
\label{shgg}
Let $G$ be the fundamental group of a finite graph of groups with
virtually polycyclic vertex groups and finite edge groups. Then
$G$ is strongly Howson and there
exists some constant $M > 0$ such that:
$$(n_1,n_2)\xi_G \leq M(n_1-1)(n_1-1)+M$$
for all $n_1,n_2 \geq 1$.
\et

\proof
By Lemma \ref{find}(i), $G$
has a finite index normal subgroup $F$ which is a finitary free
product of strongly polycyclic groups. By Theorem \ref{shspg}, there
exists a constant $M' > 0$ such that
$$(n_1,n_2)\xi_F \leq M'(n_1-1)(n_2-1)+M'$$
for all $n_1,n_2 \geq 1$. Let $m = [G:F]$ and $M = M'm^2$. By Theorem
\ref{boho}, we get
$$\begin{array}{lll}
(n_1,n_2)\xi_G&\leq&(m(n_1-1)+1,m(n_2-1)+1)\xi_F +m-1\\
&\leq&M'(m(n_1-1))(m(n_2-1))+M'+m-1\\
&=&M'm^2(n_1-1)(n_1-1)+M'+m-1\\
&\leq&M(n_1-1)(n_1-1)+M
\end{array}$$
and we are done.
\qed

\section{Takahasi's Theorem}
\label{st}

We recall, from the introduction, that a group $G$ is a {\em
Takahasi group} if every ascending chain
$$H_1 \leq H_2 \leq H_3 \leq \ldots$$
of finitely generated subgroups of $G$ with bounded rank is stationary.

Clearly, every subgroup of a Takahasi group is itself a Takahasi
group. We can prove the following partial converse:

\bt
\label{tvf}
Every finite extension of a Takahasi group is a Takahasi group.
\et

\proof
Let $G$ have a Takahasi subgroup $F$ of index $m$ and let
$H_1 \leq H_2 \leq \ldots$ be an ascending chain of subgroups of $G$
with $\rk(H_j) \leq r$ for every $j \geq 1$.

We may assume that $G$ has a standard decomposition
(\ref{stdec}). Write $K_j = H_j \cap F$. By Proposition \ref{sdsub}, there
exist $I_j \subseteq \{ 2, \ldots,m\}$ and $x_i^{(j)} \in F$ $(i \in I_j)$ such
that
\beq
\label{tvf1}
H_j = K_j \cup (\bigcup_{i \in I_j} K_jx_i^{(j)}b_i).
\eeq
Hence we have an ascending chain $K_1 \leq K_2 \leq \ldots$ of
subgroups of $F$. By (\ref{tvf1}), we have $[H_j :K_j] \leq m$ for
every $j \geq 1$. Since $\rk(H_j) \leq r$, it follows from Proposition
\ref{newrankfi} that
$$\rk(K_j) \leq m(r-1)+1.$$
Since $F$ is a Takahasi group, there exists some
$p \in \N$ such that $K_p = K_{p+1} = \ldots$.

On the other hand, we have necessarily
$$I_1 \subseteq I_2 \subseteq
\ldots \subseteq \{ 2, \ldots, m \},$$
hence there exists some $q \geq p$ and some $I \subseteq \{ 2, \ldots,
m \}$ such that
$$H_j = K_p \cup (\bigcup_{i \in I} K_px_i^{(j)}b_i)$$
for every $j \geq q$. Moreover, for every $i \in I$, we have
$$K_px_i^{(q)} \subseteq K_px_i^{(q+1)} \subseteq \ldots$$
Since two right cosets $K_px,K_py$ must be disjoint or equal, we get
$K_px_i^{(q)} = K_px_i^{(q+1)} = \ldots$ and so $H_q = H_{q+1} = \ldots$
\qed

In view of Theorem \ref{taka}, we immediately get:

\bc
\label{vftag}
Every virtually free group is a Takahasi group.
\ec

We note that, if we fix $H_1$, the length of a chain $H_1 \leq
H_2 \leq \ldots$ with subgroups of equal rank cannot be bounded, even
in the free group case:

\be
\label{cefr}
Let $A = \{ a,b,c,d,e \}$ and let $F$ be the free group on $A$. Let
$$H_1 = \langle acb\inv, ac\inv b\inv, adb\inv, ad\inv b\inv \rangle.$$
Fix $n \geq 2$ and define
$$H_i = \langle H_1, ab\inv, ae^{2^{n-i}}b\inv \rangle$$
for $i = 2, \ldots, n$. Then
$$H_1 < H_2 < \ldots < H_n$$
and all subgroups $H_i$ have rank 4.
\ee

Indeed, the Stallings automata of the $H_i$ are
of the form
$$
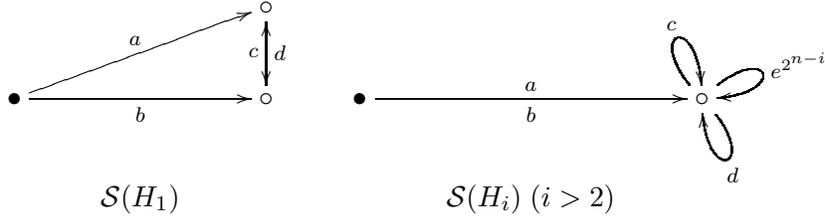
\xymatrix{
&&\circ \ar@{<->}[d]_c^d &&&\\
\bullet \ar[rr]_b \ar[urr]^a && \circ & \bullet \ar[rr]_b^a && \circ
\ar@(ur,r)^{e^{2^{n-i}}} \ar@(ul,u)^c \ar@(dr,d)^d \\
& \S(H_1) &&& {\S(H_i) \; (i > 2)} &
}$$
It follows easily from (\ref{gwp}) that $H_1 < H_2 < \ldots < H_n$. On
the other hand, (\ref{rsa}) yields
$\rk(H_1) = 6-3+1$ and
$$\rk(H_i) = (4+2^{n-i})-(1+2^{n-i})+1 = 4$$
for $i = 2,\ldots,n$.

\medskip

In the free group case, bounds can be obtained in relation with concepts such as
{\em fringe}, {\em overgroup} or {\em algebraic extension} (see
\cite{MVW}), but it is not clear how they could be efficiently generalized to
more general classes of groups.

We present now another application of Theorem \ref{tvf} which
generalizes Corollary \ref{vftag}:

\bc
\label{vpft}
The fundamental group of a finite graph of groups with
virtually polycyclic vertex groups and finite edge groups is a
Takahasi group.
\ec

\proof
Let $G$ be such a group. By Lemma \ref{find}(i), there exists a finite
index $N \unlhd G$ which is a free product of strongly
polycyclic groups. Since $G$ is finitely generated, it follows from
Proposition \ref{newrankfi} that $N$ is also finitely generated. By
Grushko Theorem, we may
write $N = S_1 \ast \ldots \ast S_t$ for some strongly polycyclic groups
$S_1, \ldots, S_t$.


Since every subgroup of a polycyclic group is finitely generated by
(\ref{shpol2}), it follows from \cite[Corollary 6.3]{Syk3} that every
ascending chain of subgroups of bounded Kurosh rank of a free
product of polycyclic groups is stationary.

Now every ascending chain of subgroups
of bounded rank of $N = S_1 \ast \ldots \ast S_t$ has also bounded
Kurosh rank by (\ref{utah}) and is therefore stationary. Thus $N$ is a
Takahasi group. By Theorem \ref{tvf}, $G$ is also a
Takahasi group.
\qed

\section{Periodic subgroups}
\label{spp}

In this section we combine Theorem \ref{tvf} with theorems on fixed
subgroups to get results on the periodic subgroups.

Given a group $G$, we denote by
$\endo(G)$ the endomorphism monoid of
$G$. Given $\p \in \endo(G)$, the {\em fixed subgroup} of $\p$ is
defined by
$$\fix(\p) = \{ x \in G \mid x\p = x\}$$
and the {\em periodic subgroup} of $\p$ is
defined by
$$\per(\p) = \bigcup_{n \geq 1} \fix(\p^n).$$
Given $x \in \per(\p)$, the {\em period} of $x$ is the least $n \geq
1$ such that $x\p^n = x$.

\bt
\label{ohio}
Let $G$ be the fundamental group of a finite graph of groups with
finitely generated virtually nilpotent vertex groups and finite edge
groups. Then there exists a constant $M > 0$ such that
$${\rm rk}({\rm Per}(\p)) \leq M$$
for every $\p \in {\rm End}(G)$.
\et

\proof
By Lemma \ref{find}(ii), $G$ has a finite index normal subgroup $N$ which is a
  finitary free product of finitely generated nilpotent groups, say $N
  = K_1 \ast \ldots \ast K_t$. Let $n
  = [G:N]$. By \cite[Lemma 2.2]{HW}, the intersection $F$ of all
  subgroups of $G$ of index $\leq n$ is a fully invariant subgroup of
  $G$, in the sense that $F\p \subseteq F$ for every $\p
  \in \endo(G)$. Moreover, since $G$ is finitely generated, we have
  $[G:F] < \infty$. Since $F \leq N$, it follows from the Kurosh
  subgroup theorem that
$F$ is isomorphic to a free product of the form
$$(\ast_{j \in J} H_j) \ast F_A,$$
where each $H_j$ is the intersection of $F$ with some conjugate of some
$K_i$. Since $G$ is finitely generated, it follows from Proposition
\ref{newrankfi} that $F$ is finitely generated and
so has finite Kurosh rank by (\ref{utah}). Similarly to the proof of Lemma
\ref{find}(ii), it follows easily that $F$ is a finitary free product
of finitely generated nilpotent groups, say $F = L_1 \ast \ldots \ast L_s$. By
\cite[Theorem 7]{Syk}, we have
$$\krk(\fix(\psi)) \leq s$$
for every $\psi \in \endo(F)$. Since each $L_i$ is polycyclic, it
follows from (\ref{shpol2}) that there exists some constant $M' > 0$
such that
$$\rk(P) \leq M'$$
for all $i \in \{ 1,\ldots, s\}$ and $P \leq L_i$. Hence we may apply
(\ref{bsvp1}) to get
$$\rk(\fix(\psi)) \leq M'\krk(\fix(\psi)) \leq M's$$
for every $\psi \in \endo(F)$.

Write $M = M's + [G:F] -1$. Let $\p \in \endo(G)$ and let $\psi =
\p|_F$. Since $F$ is a fully invariant subgroup of $G$, we have
$\psi \in \endo(F)$. Moreover, $\fix(\p) \cap F = \fix(\psi)$. By
Proposition \ref{sdsub}, we get
$$[\fix(\p):\fix(\psi)] = [\fix(\p):\fix(\p) \cap F] \leq [G:F],$$
hence \beq \label{ohio2} \rk(\fix(\p)) \leq \rk(\fix(\psi)) + [G:F]
-1 \leq M's + [G:F] -1 = M. \eeq We note that \beq \label{psv1} m|m'
\hspace{.3cm} \Rw \hspace{.3cm} \fix(\p^{m}) \leq \fix(\p^{m'}) \eeq
for all $m,m' \geq 1$: Indeed, if $m' = mk$ and $u \in
\fix(\p^{m})$, then
$$u\p^{m'} = u\p^{mk} = u\p^m\p^{m(k-1)} = u\p^{m(k-1)} = \ldots = u\p^m
= u$$
and so $u \in \fix(\p^{m'})$.

Hence we have an ascending chain of subgroups of $G$ of the form
$$\fix(\p) \leq \fix(\p^{2!}) \leq \fix(\p^{3!}) \leq \ldots$$
By (\ref{ohio2}), we have $\rk(\fix(\p^{m!})) \leq M$ for every $m
\geq 1$. Since every finitely generated nilpotent group is
polycyclic, $G$ is a Takahasi group by Corollary \ref{vpft} and so
there exists some $k \geq 1$ such that $\fix(\p^{m!}) =
\fix(\p^{k!})$ for every $m \geq k$. In view of (\ref{psv1}), we get
$$\per(\p) = \bigcup_{m \geq 1} \fix(\p^{m}) = \bigcup_{m \geq 1}
\fix(\p^{m!}) = \fix(\p^{k!}).$$
Therefore $\rk(\per(\p)) = \rk(\fix(\p^{k!})) \leq M$ by (\ref{ohio2}).
\qed

\bc
\label{boup}
Let $G$ be the fundamental group of a finite graph of groups with
finitely generated virtually nilpotent vertex groups and finite edge
groups. Let $\p \in {\rm End}(G)$. Then there exists a constant
$R_{\p} > 0$ such that every $x \in {\rm Per}(\p)$ has period $\leq R_{\p}$.
\ec

\proof
By Theorem \ref{ohio}, we have $\rk(\per(\p)) < \infty$. Assume that
$\per(\p) = \langle x_1, \ldots, x_r \rangle$. Let $R_{\p}$ denote the
least common multiple of the periods of the elements $x_1, \ldots,
x_r$. Let $x \in \per(\p)$. Then there exist $i_1, \ldots, i_n \in \{
1, \ldots, r\}$ and $\varepsilon_1, \ldots, \varepsilon_n \in \{ -1,
1\}$ such that $x = x_{i_1}^{\varepsilon_1}\ldots
x_{i_n}^{\varepsilon_n}$. It follows that
$$\begin{array}{lll}
x\p^{R_{\p}}&=&(x_{i_1}^{\varepsilon_1}\ldots
x_{i_n}^{\varepsilon_n})\p^{R_{\p}} = (x_{i_1}\p^{R_{\p}})^{\varepsilon_1}\ldots
(x_{i_n}\p^{R_{\p}})^{\varepsilon_n}\\
&=&x_{i_1}^{\varepsilon_1}\ldots x_{i_n}^{\varepsilon_n} = x,
\end{array}$$
hence $x$ has period $\leq R_{\p}$. \qed

Note that, in particular, the preceding results hold for finitely
generated virtually free groups.

We remark also that we cannot get any analogue of Theorem \ref{ohio}
involving direct products. In fact, by \cite[Theorem 4.1]{RSS}, there
exist automorphisms $\p$ of $F_2 \times \Z$ such that neither
$\fix(\p)$ nor $\per(\p)$ is finitely
generated.

\section*{Acknowledgements}

The first author is partially supported by CNPq, PRONEX-Dyn.Syst. and
FAPESB (Brazil). 

The second author acknowledges support from:
\bi
\item
CNPq (Brazil) through a BJT-A grant (process 313768/2013-7);
\item
the European Regional
Development Fund through the programme COMPETE
and the Portuguese Government through FCT (Funda\c c\~ao para a Ci\^encia e a
Tecnologia) under the project\linebreak
PEst-C/MAT/UI0144/2011.
\ei

The third author acknowledges support from a GSRT/Greece excellence
grant, cofunded by the ESF/EU and National Resources.

\end{document}